
%
%
\input amstex
\define\PP{   P_t}
\define\MD{\Cal M_k^{AH}}
\define\MG{\Cal M_k}
\define\1{1\kern-.23em\text{\rm l}}
\define\2{1\kern-.18em\text{\rm l}}
\def \qed {\hbox{\hskip5 pt\vbox{\hrule\hbox{\vrule\kern3pt\vbox{
            \kern3pt \kern3pt}\kern3pt\vrule}\hrule}}\hskip0 pt}
\magnification=1200
\baselineskip 20pt
\NoBlackBoxes

\def\sub#1 {$\underline{\hbox{#1}}$}

\documentstyle{amsppt}
\nologo
\topmatter

\title
The motive of some moduli spaces of vector bundles
over a curve.
\endtitle

\rightheadtext {The motive of some moduli spaces.}

\author
Sebastian del Ba\~no Rollin.
\endauthor

\abstract
We study the motive of the moduli spaces of semistable rank two vector
bundles over an algebraic curve.
When the degree is odd the moduli space is a smooth projective
variety,
we obtain the absolute Hodge motive of this, and in particular the
Hodge-Poincar\'e polynomial.
When the degree is even the moduli space is a singular projective
variety, we compute pure Euler characteristics and
show that only two weights can occur in each cohomology group, we also
see that its cohomology is pure up to a certain degree.
As a by-product we obtain the isogeny type of some intermediate
jacobians of the moduli spaces.
\endabstract

\date
31 January 1995
\enddate

\address {
\newline
Departament de Matem…tica Aplicada I \newline
Universitat Polit\`ecnica de
Catalunya \newline
Diagonal, 647 \newline
08028 Barcelona \newline
SPAIN \newline
e-mail:\tt Sebastian\@ma1.upc.es }
\endaddress

\endtopmatter

\document

\footnote""{Partially supported by DGCYT grant PB93-0790.\hfill }

\heading \S 0.Introduction.\endheading

The moduli space of stable vector bundles over an algebraic curve
is a relatively well-known object, it has received great attention for
the last twenty years, in particular when the rank and degree are
coprime its cohomology has been shown to be
torsion free and its Betti numbers are known. However the methods used
in studying its cohomology are
topological (\cite {Ne}),
number theoretical
(\cite {HN}) or
infinite-dimensional (\cite {AB}, see also
\cite {LPV}),
and these, at least in
principle, do not yield any information of the analytic/algebraic
structure of the moduli space and the dependence of this on the curve,
for instance its Hodge
numbers remain unknown.

In this paper we use a
recent construction by M.Thaddeus (\cite {Th}) to give
a description of the motivic Poincar\'e polynomial of the moduli space
of
rank two semistable vector bundles of fixed determinant on an algebraic
curve.
It is an idea of Grothendieck (see \cite {S}) that one should work in
the Grothendieck group  $K_0$ of the category of motives, this is
where the motivic Poincar\'e polynomial lives.
We believe that the theory of
motives is an effective language to express clearly and precisely how
the algebro-geometric properties  of  $C$
influence those of the moduli space  $N_0(r,d)$.
However at the present moment we do not have at our disposal
the true category of motives  $\MG $ of Grothendieck, since the Standard
Conjectures remain unproven,
so we use the definition by Deligne of absolute Hodge motives
$\MD$.

We start by giving a quick review of the theory of absolute
Hodge motives, the natural
language in which motives are expressed is that of tannakian categories
so we
recall the basic facts of these,
we
also define the motivic Poincar\'e polynomial,
this is done in \S 1.

In order to carry out the calculations we need a motivic version of
MacDonald's formula for the Betti numbers of a symmetric power, we
do this in \S2. Note that in fact we get an expression for the motive of
$X^{(n)}$ and not only  $\PP X^{(n)}\in K_0\MD $.

Then in \S 3 we give a short account of
Thaddeus' construction of the moduli spaces of pairs.

In \S 4 is where with the aid of Thaddeus' construction we manage to
calculate the motivic Poincar\'e polynomial of the moduli space
$N_0(2,1)$ of stable rank  $2$ vector bundles with fixed odd
determinant.

In \S 5 we study the singular moduli space  $N_0(2,0)$ of rank  $2$
semistable vector bundles with fixed even determinant. Our results difer
from
those of
Kirwan
(\cite {K})
in that we use pure
Poincar\'e polynomials
whereas she  works
with a canonical partial desingularization of the moduli space to get
the intersection cohomology Poincar\'e polynomial.

Finally in \S 6 we extract information concerning the intermediate
jacobians of the moduli spaces from the motivic Poincar\'e polynomial.

\demo{Acknowledgement}
I would like to express my gratitude to Vicente Navarro and Pere Pascual
for the invaluable help and encouragement they have shown me during
the preparation of this paper. I would also like to thank my colleagues
at the Barcelona Cigrons Geometry Seminar.
\enddemo

\bigbreak
\heading
\S 1. Absolute Hodge motives.
\endheading

In this section we give a brief review of
the theory of absolute Hodge motives and related topics. For proofs of
theorems and more precise statements refer to \cite {DM}.

Let  $k$ be a field of characteristic zero.

\subheading
{Tannakian Categories}

By a tannakian category we shall mean a  $k$-linear abelian neutral
rigid
tensor category with  $End(\1)=k$.
Let $\Cal C$ be a tannakian category,
the fact that $\Cal C$ is neutral means
there is a faithful exact
functor  from $\Cal C$ to the category of finite dimensional vector
spaces over  $k$

\centerline {$\omega :\Cal C\longrightarrow \bold {Vec_k} $}

\noindent
called a fibre functor.

In a rigid tensor category there is a concept of rank. However in the
case of a tannakian category this definition of rank gets simplified by
the use
of the fibre functor, $\text{rank}M=\text{dim}_k\omega (M)$.

An example of tannakian category is the category of
finite dimensional $k$-re\-pre\-sen\-tations
of an affine group scheme  $\Cal G$ over  $k$, the fibre functor being
the obvious one.
In fact a fundamental theorem
states that all tannakian categories arise in this way, so that if
$\Cal C$ is a tannakian category with fibre functor  $\omega $ then
there exists an equivalence of categories  $\Cal C\longrightarrow \bold
{Rep\text{-}\Cal G}$ compatible with the fibre functors,  $\Cal G$ is
then called the Galois group or fundamental group of the tannakian
category  $\Cal C$, $Gal(\Cal C)$.

\demo {Examples}

\item{\it 1.}
Consider the category of finite dimensional graded vector spaces over
$k$, it
is a tannakian category which is easily seen to be equivalent to  $\bold
{Rep\text{-}\Bbb G_m}$.

\item{\it 2.}
The category of local systems of finite dimensional complex vector
spaces over a topological space  $X$
is a tannakian category. A fibre functor is obtained by assigning to
each local system  $\Cal L$  the complex vector space  $\Cal L_x$, where
$x$ is a point of  $X$.

\item{\it 3.}
The category of rational pure Hodge structures, $\Cal {HS}_{\Bbb Q}$ is
a tannakian category.
\enddemo

There is also the notion of a graded tannakian category. It is one where
every object has a direct sum decomposition compatible with Knneth.
This is better expressed by using the Galois group, a graded
tannakian category is a tannakian category  $\Cal C$ together with
a central morphism  $\Bbb G_m\longrightarrow Gal(\Cal C)$.

However a richer structure appears quite naturally in the theory of
motives:
$(\Cal C,w,T)$ is called a  Tate triple if   $\Cal C$ is a tannakian
category graded by $w:\Bbb G_m\longrightarrow Gal(\Cal C)$, and  $T$ is
a weight  $-2$ invertible object called the Tate object. The
result of tensoring
an object by the Tate object is usually refered to as a Tate twist.
A standard way to abreviate
$A\otimes T^{\otimes i}$ is  $A(i)$.

The following definition will be useful

\proclaim{Definition} If  $\Cal C$ is a tensor category then  $\Cal
C[[T]]$ is the tensor category whose objects are
$$Ob(\Cal C[[T]])=\{ (A_i)_{i\in\Bbb N}\vert A_i\in Ob(\Cal C)\}$$
often written  as $\sum A_iT^i$. The
morphisms are defined by
$$Hom (\sum A_iT^i,\sum B_iT^i)=\prod Hom(A_i,B_i) $$
there is a natural functor
$$\underset T^n\to {\text{Coef}}:\Cal C[[T]]\longrightarrow \Cal C$$
sending  $\sum A_iT^i$ to  $A_n$, it is not however a tensor functor.
\endproclaim

Recall that in a tensor category there are commutation constraints, that
is for every pair of objects $M,N\in Ob(\Cal C)$ isomorphisms
$$
M\otimes N \buildrel \varphi \over \longrightarrow
N\otimes M
$$

Let  $\Cal C$ be a graded tensor category. Consider
the
new commutation constraints, given on pure degree objects  $M_i$,  $N_j$
of weights  $i$ and  $j$ by
$$\align
\psi :M_i\otimes &N_j \longrightarrow N_j\otimes M_i \\
&\psi =(-1)^{ij}\varphi
\endalign $$
where  $\varphi $ are the old commutation constraints. Call  $\buildrel
\centerdot \over {\Cal C} $ the resulting tensor category.

In the case when  $\Cal C$ is a tannakian category
 $\buildrel \centerdot \over {\Cal C}$ need not be tannakian.
For instance, if  $\Cal C$ is $\Cal M_k^{AH}$ or  $\Cal M_k$
then $\overset \centerdot \to {\Cal C}$ is
called the false category of motives  $\buildrel \centerdot \over {\Cal
M^{AH}_k}$ or  $\buildrel \centerdot \over {\Cal M}_k$.

\medbreak
\subheading {Absolute Hodge Motives}

Let  $k$ be a field of characteristic zero embeddable in  $\Bbb C$. We
shall work with the category of smooth projective varieties over  $k$,
$\Cal V_k$.

The main problem in the theory of motives is to find a
tannakian category that factors all possible cohomology functors.
Grothendieck gave a
construction of such a category  $\MG $ (see \cite {Ma}) but in order to
prove it has
the required properties one needs the Standard Conjectures which remain
unproven.

Deligne (\cite {DM}) has given a temporary working definition for
motives,
these are the absolute Hodge motives which we shall use in what follows.
The category $\MD $ is constructed in
exactly the same way as  $\MG $ but using absolute
Hodge cycles instead of algebraic cycles.
We recall that an absolute Hodge cycle of  $X$  of
codimension  $p$  is an element of
$$
F^0H^{2p}_{DR}(X)(p)\times \underset l\to \prod
H^{2p}_{\text{\'et}}(\overline X,\Bbb Q_{\ell })(p)\times
\underset \sigma:k \hookrightarrow \Bbb C \to \prod
 H^{2p}_{sing}(X_\sigma ,\Bbb Q)(p)
$$
such that it is compatible with the comparison
isomorphisms. We denote the group of such cycles by $Z^p_{AH}(X)$

In the same manner as with Grothendieck motives we get
a functor, $h$, from
the category of smooth projective varieties over  $k$ to the category of
AH-motives.

One advantage of working with AH-motives is that the Knneth components
of the diagonal in  $H^{2d}(X\times X)$ are again AH-cycles, so we get a
decomposition  $hX=h^0X\oplus h^1X\oplus \cdots \oplus h^{2d}X$. This
makes $\MD $ into
a graded tannakian category, it is customary to refer to this grading as
the weight grading. A motive that is zero in all degrees except in
one is called a pure weight motive.
As $\MD $ is a graded tannakian category one has a graded fibre functor
$$ \align
\MD   &\longrightarrow \bold {Grad\text{-}Vec}_k \\
M=\oplus M_i&\longrightarrow \oplus H^i_{DR}(M)
\endalign $$

It is proven in \cite {DM} that  $\MD $ is in fact a polarized Tate
triple, the
Tate object is $\1 (1)$. Using the fact that  $\MD $ is polarized, one
can prove that the category  $\MD $ is semisimple, that is, every motive
is the direct sum of simple motives.

\demo{Remarks}

\item {\it 1.} The cycle maps  $Z^p(X)\longrightarrow H^{2p}(X)(p)$
produce an absolute Hodge cycle for each p-codimensional algebraic
cycle so one gets a morphism  $Z^p(X)\longrightarrow Z_{AH}^p(X)$.
This way we get a functor  $\MG \longrightarrow \MD$.

\item {\it 2.}
An important thing to know about the category of
AH-motives is that it is a full subcategory of the category of
realization systems defined in \cite {J}. The motive $h^i(X)$ can thus
be seen as a triple  $\left(H_{DR}^i(X,k),
H_\sigma ^i(X,\Bbb Q),
H^i_{\text{\'et}}(\overline X,\Bbb Q_\ell)\right)$,
where $H_{DR}^i(X,k)$ is a finite dimensional
$k$-vector space with a filtration (the Hodge filtration),
for each embedding  $k\buildrel \sigma \over \hookrightarrow \Bbb C$
$H_\sigma
^i(X,\Bbb Q)$ is a rational pure Hodge structure of weight $i$ and
for each prime $\ell $  $H^i_{\text{\'et}}(X,\Bbb Q_\ell )$ is a
$Gal(\overline k,k)$-module, together with comparison isomorphisms.
\enddemo

\medbreak
\subheading
{$\bold {K_0\MD} $ and the motivic Poincar\'e polynomial}

Recall that to every abelian category $\Cal C$  one can attach the
Grothendieck group $K_0\Cal C$. Moreover if  $\Cal C$ is a (graded)
tensor
category then $K_0 \Cal C$ is a (graded) unitary commutative ring. Given
$A$ an object
$\Cal C$ we shall use the notation $[A]$ for its image in  $K_0\Cal C$.

In particular if we put  $\Cal C=\MD $ we get a graded ring $K_0\MD $.
In the category  $\MD $ one has the Tate twist
$$\align
\MD &\buildrel \cdot \otimes \2(n)\over \longrightarrow \MD \\
A&\longrightarrow A(n)
\endalign $$
and the dualising functor
$$\align
\MD &\buildrel \cdot ^{\vee } \over \longrightarrow \MD \\
A&\longrightarrow A^{\vee }
\endalign $$
both of which are exact functors so they descend to additive morphisms
of the graded ring  $K_0\MD $.

\proclaim{Definition }
Let  $M$ be an AH-motive and  $M=\oplus M_i$ its weight grading, then
its motivic Poincar\'e polynomial is defined to be its class in the
graded ring  $K_0\MD $
$$
P_t^{\text{mot}} (M) =[M]=\sum [M_i]\in K_0\MD
$$
\endproclaim

Note that this is not really a polynomial, it is an element of a graded
ring.
We shall often drop the ${}^{\text{mot}}$ and just
write  $\PP $.
If  $X$ is a smooth projective variety over  $k$ then we shall write
$\PP (X)=\PP (hX)$.

This is a generalisation of the usual Poincar\'e polynomial as can be
seen by following  $h(X)=\oplus h^i(X)$ through the commutative diagram
$$ \CD
Ob( \MD) @>[\cdot ]>> K_0\MD  \\
@V\omega VV @V K_0(\omega) VV        \\
Ob(\bold {Grad\text{-}Vec_k}) @>[\cdot ]>>
K_0\bold {Grad\text{-}Vec_k }=\Bbb Z[t,t^{-1}]
\endCD
$$

In general the map  from isomorphism classes of objects of an abelian
category  $\Cal C$ to $K_0\Cal C$
is not injective, but as an application of the fact that
$\MD $ is semisimple we show now that
this is the case for  $\MD $.

The following proposition is easily proven.

\proclaim {Proposition} Let  $\Cal C$ be an artinian abelian semisimple
category (for ex\-am\-ple  $\MD $) and  $A,B,C\in Ob(\Cal C)$. If
$A\oplus C\simeq B\oplus C$ then $A\simeq B$. \endproclaim

\proclaim{Corollary} Let  $M,N$ be AH-motives. If
$\PP M=\PP N $
then  $M\simeq N$.
\endproclaim

\demo{Proof}
$\PP M=\PP N$ means that there exists  $P$
with
$M\oplus P\simeq N\oplus P$, now
use the previous proposition.
\hfill\qed
\enddemo

Therefore whenever we need to prove an equality of motives it will be
enough to prove it in  $K_0$, and this is normally easier to write.

\medbreak
\subheading {Mixed absolute Hodge motives}

The geometric methods in the definition of AH-motives do not extend at
the present moment
to the case of open or singular varieties. As already
mentioned  $\MD $ is a full subcategory of  the category of
realization systems $\Cal R_k$, this is very useful to construct a
category of mixed absolute Hodge motives as there is a reasonable
candidate for category of mixed systems of realizations,  $\Cal {MR}_k$
together with natural functors  $h^i:\Cal W_k\longrightarrow \Cal
{MR}_k$. Where  $\Cal W_k$ denotes the category of varieties over  $k$
(not necessarily smooth or proper).

Let  $\Cal V^0_k$ denote the category of smooth varieties over  $k$ (not
necessarily proper).
Jannsen (\cite {J}) defines  $\Cal {MM}_k^{AH}$ to be the full tannakian
subcategory of  $\Cal {MR}_k$ generated by the image of the  $h^i:\Cal
V^0_k\longrightarrow \Cal {MR}_k$.

There is a functor

\centerline {$h:\Cal W_k\longrightarrow \bold {Grad\text{-}\Cal
{MM}_k^{AH} }
$}

\noindent
which assigns  $\oplus h^i(X)$ to the variety  $X$

There is a natural fully faithful functor  $\MD \longrightarrow \Cal
{MM}_k^{AH}$,  $\MD $ can thus be seen as a full subcategory of  $\Cal
{MM}^{AH}_k $.
The fact that every object
in $\Cal {MM}_k^{AH}$ is an extension of objects in  $\MD $ implies that
the previous functor induces an isomorphism  of rings $K_0\MD
\buildrel \simeq \over \longrightarrow K_0\Cal {MM}_k^{AH}$.

We now define a polynomial which via the mentioned isomorphism extends
the motivic Poincar\'e polynomial.

\proclaim{Definition (\cite {Na})}
Let  $M=\oplus M_i\in Ob( \bold{Grad\text{-}}\Cal {MM}^{AH}_k)$ be a
graded mixed motive, then the pure motivic Poincar\'e polynomial is

\centerline {$
P^{\text{mot}}_t(M)=\underset m\to \sum \left(
\underset i\to \sum (-1)^{m+i}[Gr^{W_\cdot}_m M_i] \right)
\in K_0\MD
$}
\endproclaim

\demo{Remarks}

\item {\it 1.} If  $M=\oplus M_i$ with $M_i$ a pure motive of
weight i, $Gr^{W_\cdot}_mM_i$
is equal to  $M_i$ if  $i=m$ and zero if $i\neq m$ so that this
polynomial coincides with the one already defined.

\item {\it 2.} Note that in the mixed case  $\PP M$ does not coincide
with the class of  $M$ in  $K_0\MD $.

\item {\it 3.} Let  $X$ be a variety over  $k$,
$\oplus
h^iX $
its mixed motive and $P_t^{\text{mot}}$
its motivic Poincar\'e polynomial.
Composition with the ring morphism
$$
K_0\MD
\longrightarrow K_0 \bold {Grad\text{-}Vec}_{k}=\Bbb Z [t,t^{-1}]
$$
does not yield the
classical Poincar\'e polynomial  $P_tX:=\sum \text{dim} H^i(X,\Bbb
Q)t^i$ but rather the pure Poincar\'e polynomial defined by
$$
P_t^{pur}(X)=\underset m\to \sum
\chi _m^{pur}(X)t^m
\text{, where }
\chi_m^{pur}(X)=\underset i\to \sum
(-1)^{i+m}\text {dim}Gr_m^{W_\cdot }H^i(X,\Bbb Q)
$$
(c.f. \cite {G,185-191}, \cite{S} and \cite {Na})
this is better suited for computations than the ordinary Poincar\'e
polynomial. For example if  $Y$ is a closed subvariety of  $X$ of
codimension  $d$ and both  $X$ and  $Y$ are smooth one has the Gysin
exact sequence,
$$
\cdots \longrightarrow
h^{i-2d}Y(-d)\longrightarrow h^iX\longrightarrow
h^i(X-Y)\longrightarrow \cdots
$$
as the functor  $Gr_m^{W_\cdot}$ is exact one gets an equality in
$K_0\MD $
$$
\underset i\to \sum
(-1)^i
[Gr_m^{W_\cdot }h^iX]
=
\sum (-1)^i [Gr_m^{W_\cdot }h^i(X-Y)]
+
\sum (-1)^i [Gr_{m-2d}^{W_\cdot }h^{i-2d}Y](-d)
$$
so that  $P^{mot}_tX=P^{mot}_t(X-Y)+P_t^{mot}Y(-d)$.

\enddemo

\bigbreak
\heading
\S 2. A motivic MacDonald formula.
\endheading

Let  $X$ be a compact polyhedron and consider
 $X^{(n)}$ the symmetric power of  $X$, this is the quotient of  $X^n$
by the natural action of the
symmetric group $\frak S_n$.
MacDonald gave a formula (\cite {M}) that
produces Betti numbers of
$X^{(n)}$ in terms of those of $X$, explicitly
$$
P_tX^{(n)}=
{\underset T^n\to {\text{Coef}} }
{   (1+tT)^{b_1(X)}\cdot (1+t^3T)^{b_3(X)}\cdots
    \over
    (1-T)^{b_0(X)}\cdot (1-t^2T)^{b_2(X)}\cdots
}
$$

In this paragraph we intend to give a motivic version of
MacDonald's formula valid in any neutral $k$-linear graded tannakian
category, in particular that of Absolute Hodge Motives or
conjecturally Grothendieck's category of pure motives.

Throughout all this section $k$ will denote a field of characteristic
$0$.
Let  $\Cal C $ be a
tannakian category,
in \cite{DM, pg.106} it is shown that the commutation
constraints
can be extended in a natural way to cover the case of more than two
factors,
for every   $\sigma \in \frak S_n$, we get isomorphisms
$$
\varphi _\sigma :M_1\otimes\cdots \otimes M_n
\longrightarrow
M_{\sigma ^{-1}(1)}\otimes \cdots \otimes M_{\sigma ^{-1}(n)}
$$
in particular
if $M\in
Ob(\Cal C)$
this defines an action of   $\frak S_n$ on  $M^{\otimes n}$
$$
\frak S_n\buildrel \varphi \over \longrightarrow Aut(M^{\otimes n})\text
{.}
$$

Let  $\varepsilon :\frak S_n\longrightarrow \{+1,-1\}$ denote the
signature

\proclaim{Definition}
Recall that a tannakian category is abelian, in particular any morphism
has an image. Given a  $M\in Ob(\Cal C)$ define
$S^iM$ (resp.  $\wedge ^iM$) to be the image of the morphism
${1\over i!}\underset \sigma \in
\frak S_i\to \sum \varphi _\sigma :M^{\otimes i}\longrightarrow
M^{\otimes i}$ (resp.
${1\over i!}\underset \sigma \in
\frak S_i\to \sum \varepsilon (\sigma ) \cdot \varphi _\sigma $).
Extend this definition to the case  $i=0$
by putting $S^0M=\wedge ^0M=\1 $.
\endproclaim

\demo {Remarks}

\item {\it 1.}
As the fibre functor $\omega $ is a tensor functor it sends  $\varphi
_\sigma $ to the canonical commutation constraints in  $\bold
{Vec}_k$, this combined with the fact that  $\omega $ is exact
gives immediately that $\omega \left(
\wedge^iM\right)=\wedge ^i\omega (M)$.
Using the faithfulness of  $\omega $  we see
that  $\wedge ^iM=0$ for  $i>\text {rank}M$.

\item {\it 2.} If  $M$ is a rank one object then using again the fibre
functor one sees that  $S^iM=M^{\otimes i}$.
\enddemo

\proclaim{Definition}
If $M\in Ob(\Cal C)$ define
$$ \align
(1+T)^M&=\sum \wedge ^iM\cdot T^i\in Ob(\Cal C[T]) \\
(1-T)^{-M}&=\sum S ^iM\cdot T^i\in Ob(\Cal C[[T]])
\endalign $$

If the rank of  $M$ is one then  $(1-T)^{-M}=\sum M^{\otimes i}T^i$ so
we shall also use the notation  ${1\over 1-MT}$ in this case.
If  $A$ is an invertible object then
${1\over A-BT}$ will mean the reasonable thing
$${1\over A-BT}=A^{-1}\otimes \left( {1\over 1-A^{-1}BT }\right)\in
Ob(\Cal C[[T]]) $$
\endproclaim

For the rest of the section  $\Cal C$ will denote a graded tannakian
category over  $k$. Recall that for a graded tensor category  $\Cal C$
we
defined in \S 1 a tensor category  $\buildrel \centerdot \over {\Cal C}$
by changing certain signs in the commutation constraints.

\proclaim{Definition} Define the symmetric power of  $M$,
$M^{(i)}$,
the same way as  $S^iM$
but using the commutation
constraints from  $\overset \centerdot \to {\Cal C}$.
\endproclaim

\proclaim{Proposition} Let  $M$ be a pure degree object of weight  $n$
then
 $M^{(i)}$ is  $S^iM$ if  $n$ is even and  $\wedge ^iM$ if  $n$ is odd.
\endproclaim

\demo {Proof}
If  $M$ is pure of even weight then the commutative constraints
$$ M\otimes M\longrightarrow M\otimes M
$$
are the same in  $\buildrel \centerdot \over {\Cal C}$ and in  $\Cal C$
so  $M^{(n)}=S^nM$.

In the odd weight case the commutation constraints change sign and when
more than one factor appears then the sign is given by the signature
$\varepsilon $ so we get  $M^{(n)}=\wedge ^nM$.
\hfill \qed \enddemo

The next theorem gives an expression for the symmetric power of an
object in terms of symmetric powers of its pure components, it is our
motivic version of the MacDonald formula.

\proclaim{Theorem} Let  $M=\oplus M_i$ be an object in a graded neutral
$k$-linear tannakian category, then
$$
M^{(n)}=
\underset T^n \to {\text{Coef}}
{
 \cdots
\otimes
(1+T)^{M_{-1}}
\otimes (1+T)^{M_1}
\otimes (1+T)^{M_3}\otimes \cdots
\over
\cdots
\otimes
(1-T)^{M_{-2}}\otimes
(1-T)^{M_0}\otimes
(1-T)^{M_2}\otimes \cdots
}
$$
\endproclaim

\demo{Proof}
We need to see that
$M^{(n)}=(M^{\otimes n})^{\frak S_n}$
is
$$
\align
\underset T^n \to {\text{Coef}}
&
\left(
\cdots \otimes
\underset i \to \sum \wedge ^iM_{-1} T^i\otimes
\underset i \to \sum \wedge ^iM_1 T^i\otimes
\underset i \to \sum \wedge ^iM_3 T^i\otimes \cdots
\right.
\\
&
\left.
\otimes
\underset i \to \sum S^iM_{-2} T^i\otimes
\underset i \to \sum S^iM_{0} T^i\otimes
\underset i \to \sum S^iM_2 T^i\otimes \cdots
\right)
\\
&=\underset \lambda _1+\cdots +\lambda _k=n \to \sum
M_{r_1}^{(\lambda _1)}\otimes \cdots \otimes M_ {r_k}^{(\lambda _k)}
\endalign
$$

On the other hand, by Knneth
$$
\align
M^{\otimes n}
&=
\underset r\in \Bbb Z \to \sum
\left(
\underset {}_{r_1+\cdots +r_n=r} \to \bigoplus
M_{r_1}\otimes \cdots \otimes M_{r_n}
\right)
\\
(M^{\otimes n})_r
&=
\underset {}_{r_1+\cdots +r_n=r} \to \bigoplus
M_{r_1}\otimes \cdots \otimes M_{r_n}
\endalign
$$
using next lemma we get that the  $\frak S_n$-invariant subobject is
isomorphic to
$$
\underset
{
\underset {s_1<\cdots <s_k} \to
{  {\sum \lambda _is_i=r}
\atop {\sum \lambda _i=n}  }
}
\to
\bigoplus
\left(
M_{s_1}\otimes
\buildrel \lambda _1\over \cdots
\otimes M_{s_1}\otimes
\cdots
\otimes M_{s_k}\otimes
\buildrel \lambda _k\over \cdots
\otimes M_{s_k}
\right)
^{
\frak S_{\lambda _1}
\times \cdots
\times \frak S_{\lambda _k}
}
$$

Now according to whether $s_i$ is even or odd, the action of  $\frak
S_{\lambda
_i}\hookrightarrow \frak S_n$ on $M_{s_i}\otimes \cdots \otimes
M_{s_i}$ is the canonical
($\varphi $) or the anticanonical one ($\varphi
\cdot \varepsilon $) and thus  $(M_{s_i}^{\otimes \lambda _i})^{\frak
S_{\lambda _i}}$ is  $S^{\lambda _i}M_{s_i}$ or either  $\wedge
^{\lambda _i}M_{s_i}$. This proves the theorem.
\hfill \qed
\enddemo

\proclaim{Lemma} There is a natural isomorphism
$$
(M^{\otimes n})_r^{\frak S_n}
\simeq
\underset
{
\underset {r_1<\cdots <r_k} \to
{  {\sum \lambda _ir_i=r}
\atop {\sum \lambda _i=n}  }
}
\to
\bigoplus
\left(
M_{r_1}\otimes
\buildrel \lambda _1\over \cdots
\otimes M_{r_1}\otimes
\cdots
\otimes M_{r_k}\otimes
\buildrel \lambda _k\over \cdots
\otimes M_{r_k}
\right)
^{
\frak S_{\lambda _1}
\times \cdots
\times \frak S_{\lambda _k}
}
$$
\endproclaim

\demo{Proof}
The inclusion  $(M^{\otimes n})_r^{\frak S_n}\hookrightarrow
(M^{\otimes n})_r$ followed by the obvious projection gives a morphism

\centerline {$
(M^{\otimes n})_r^{\frak S_n}
\longrightarrow
\underset
{
\underset {r_1<\cdots <r_k} \to
{  {\sum \lambda _ir_i=r}
\atop {\sum \lambda _i=n}  }
}
\to
\bigoplus
\left(
M_{r_1}\otimes
\buildrel \lambda _1\over \cdots
\otimes M_{r_1}\otimes
\cdots
\otimes M_{r_k}\otimes
\buildrel \lambda _k\over \cdots
\otimes M_{r_k}
\right)
$}

Applying the exactness and faithfulness of the fibre functor one sees
that the image is contained in
$$
\underset
{
\underset {r_1<\cdots <r_k} \to
{  {\sum \lambda _ir_i=r}
\atop {\sum \lambda _i=n}  }
}
\to
\bigoplus
\left(
M_{r_1}\otimes
\buildrel \lambda _1\over \cdots
\otimes M_{r_1}\otimes
\cdots
\otimes M_{r_k}\otimes
\buildrel \lambda _k\over \cdots
\otimes M_{r_k}
\right)
^{
\frak S_{\lambda _1}
\times \cdots
\times \frak S_{\lambda _k}
}
$$
and that in fact the induced morphism from  $(M^{\otimes n})_r^{\frak
S_n}$
to this is injective. So it is enough to see that they both have the
same rank but this is just the statement of MacDonald's main theorem
(\cite {M}). \hfill \qed
\enddemo

Let  $X$ be a smooth projective variety over  $k$.
Write
 $h(X)=\oplus h^i(X)\in Ob(\MD )$,
then by
Proposition 6.8 in \cite {DM}
$h(X^{(n)})=(h(X)^{\otimes n})^{\frak S_n}$, where the action of  $\frak
S_n $ is the one arising from the geometric commutations
$X\times \cdots \times X\longrightarrow X\times \cdots \times X$,
so that
$(h(X)^{\otimes n})^{\frak S_n} =
h(X)^{(n)}$
and this is computed using the formula in the Theorem, so we have

\proclaim{Corollary}
$$
h(X^{(n)})=\underset T^n\to {\text{Coef}}
{
(1+T)^{h^1X}\otimes
(1+T)^{h^3X}\otimes \cdots
\over
(1-T)^{h^0X}\otimes
(1-T)^{h^2X}\otimes \cdots
}
$$
and if  $C$ is a smooth projective curve
$$
h(C^{(n)})=\underset T^n\to {\text {Coef}}
{
(1+T)^{h^1C}
\over (1-\1T)(1-\1(-1)T)
}
$$
\endproclaim

\demo{Remark} If we apply the graded fibre functor
$$
\MD
\buildrel
H_{DR}^* \over
\longrightarrow \bold { Grad \text {-}Vec_{k}  }
$$
followed by
$
Ob(\bold { Grad \text {-}Vec_{k}  })
\buildrel [\cdot ]\over \longrightarrow \Bbb Z[t,t^{-1}]$
we get the  classical MacDonald formula, whereas if we do the same
with
$$
\MD
\buildrel H^{*,*}_{DR} \over
\longrightarrow \bold { BiGrad \text {-}Vec_{k}  }
$$
and
$
Ob(\bold { BiGrad \text {-}Vec_{k}  })
\buildrel [\cdot ]\over \longrightarrow \Bbb Z[x,y,x^{-1},y^{-1}]$
we get the formula for the Hodge numbers in \cite {Bu}.
\enddemo

\bigbreak
\heading
\S 3. Thaddeus' construction.
\endheading

In this section we review the basic construction of Thaddeus we shall
use, for a more complete exposition see \cite{Th}.

Let  $C$ be a fixed smooth projective algebraic curve of genus  $g\geq
2$ over  $k$ an algebraically closed field of zero characteristic and
$\Cal L$
a line bundle over C of {\it large} degree $d$. The moduli spaces we are
primarily interested in are  $N_0(2,d)(C)$ the moduli space of rank
$2$  semistable vector bundles with fixed determinant over  $C$.
They depend on the curve  $C$ however we shall simply write  $N_0(2,d)$.

Thaddeus
considers the problem
of giving a moduli space for pairs  $(E,s)$, where $E$ is a rank  $2$
vector bundle over the curve
$C$ with fixed determinant  $\Cal L$  and  $s$ is a non-zero section of
$E$. It appears that there are many possible
definitions for stability of a pair depending on a parameter
$\sigma \in [0,{d\over 2}]$, for  $\sigma $ varying in certain
open disjoint intervals there are no strictly semistable pairs and
one obtains a finite list of fine moduli spaces  of pairs $M_0$, ...,
$M_\omega $ ($\omega =\left[ {d-1\over 2}\right]$).

These different moduli spaces are all birrational and are related by
a special kind of birrational maps called flips. In
this context a flip between two varieties  $X$ and  $Y$ means that  $X$
and  $Y$ have a common blow-up,  $\tilde X\simeq \tilde Y$, with the
same exceptional locus.
Luckily the centers of these blow-ups are nonsingular so that one can
use the standard formula for the Poincar\'e polynomial of a blow-up.
Of course, to be able to work out Poincar\'e polynomials of the moduli
spaces we need to know the centers of these blow-ups, these turn out
to be  a couple of subvarieties of $M_j$ called $\Bbb PW_i^+$ and
$\Bbb PW_{i+1}^-$ isomorphic to  certain projective bundles over
symmetric products of the curve:
$\Bbb PW_i^+$ is a
$\Bbb P^{d-2i+g-2}$-bundle over  $C^{(i)}$
and $\Bbb PW_{i+1}^-$ is a
$\Bbb P^{i}$-bundle over  $C^{(i+1)}$.
To summarize,  the blow-up of $M_{i}$ along  $\Bbb PW_i^-$ is
isomorphic to the blow-up of  $M_{i-1}$ along  $\Bbb PW_i^+$.
We can picture this chain of flips:

$$\matrix
&&\tilde M_2 &&&& \tilde M_3   &&& \tilde M_\omega  \\
&\swarrow &&\searrow && \swarrow &&\searrow &\cdots \swarrow &&\searrow
\\
M_1 &&&& M_2&&&&&&& M_\omega  \\
\downarrow &&&&&&&&&&& \downarrow \\
M_0  &&&&&&&&&&&  N_0(2,d)\\
\endmatrix $$

Moreover, it is easy to see
that $M_0$ is a projective
space of dimension  $d+g-2$. In the other extreme we have $M_\omega $,
in the case
when $\text{deg}\Cal L$ is odd $M_\omega $ is a projective bundle of
relative dimension  $d-2g+1$ over $N_0(2,1)$,
the moduli space of rank two stable vector bundles over  $C$ with fixed
odd determinant, whereas if  $\text {deg}\Cal L$ is even we have a map
from
$M_\omega $ to the analogous moduli space which is only a projective
fibration over the stable locus.

The above description of the centers of the flips enables Thaddeus to
recover the formula by Harder-Narasimhan for the Poincar\'e polynomial
of $N_0(2,1)$.
In the next sections we shall use the formulae for the motive of
blow-ups and our motivic MacDonald formula in the same way to get
an expresion for the motive of  $N_0(2,1)$ and  $N_0(2,0)$.

\bigbreak
\define \SSS {\underset j=0\to {\buildrel i\over \sum }}
\heading
\S 4. The motive of $N_0(2,1)$.
\endheading

The purpose of this section is to give an expression for  $hN_0(2,1)$ in
terms of  $h^1C$ and  $\1(1)$. An immediate consequence is that
$hN_0(2,1)$ is in the tannakian subcategory of  $\MD$
generated by  $hC$ and  $\1(1)$.

The calculation of the Poincar\'e polynomial of the moduli
space
involves some infinite sums of motives thus falling
outside of the ring $K_0\MD $, to formalise this we need
to construct a greater ring $\widehat  {K_0\MD }$, the ring of Laurent
series of motives, this is done as follows: first consider the subring
$$
K_0{\MD } ^+=\{ x\in K_0\MD \vert \text {deg}x\geq 0\}\subset K_0\MD,
$$
complete it with respect to the ideal  $I$  formed by the strictly
positive
degree elements, tensor the result by  $K_0 \MD$ over  $K_0 {\MD } ^+$,
then the result is the ring we were looking for,
$$
\widehat  {K_0 \MD }=K_0\MD \underset K_0{\MD }^+\to \otimes
\widehat  {K_0{\MD}} ^+_I.
$$

If $A,B$ are invertible motives, with  $\text{deg}B>\text{deg}A$, and
$C$
is any motive then define  ${C\over
A-B}$ to be the reasonable thing in $\widehat  {K_0\MD }$,
that is
$$ {C\over A-B}:=C\cdot (A^{-1}+ A^{-2}\cdot B+ A^{-3}\cdot
B^2+
\cdots ).$$

\proclaim{Proposition}
The motive of  the moduli space of pairs $M_i$ is given by
$$
hM_i=
\underset j=0 \to {\buildrel i\over \sum }
hC^{(j)}\otimes \left(
\1(-j)\oplus \cdots \oplus \1(-d+2j-g+2)
\right)
$$
and its motivic Poincar\'e polynomial is
$$ \align
\PP ^{\text {mot}}M_i=
&{\1 \over \1-\1 (-1)}\cdot
\underset T^{i}\to {\text {Coef}}
\left(
{\1 (-d+2i-g+1)\over \1(-2)T-\1}-
{\1 (-i-1)\over T-\1 (-1)}
\right)  \cdot
\\
& \cdot
{(1+T)^{h^1C}\over (\1-T)(\1-\1(-1)T)}.
\endalign $$
\endproclaim

\demo{Proof}
Notice that we suppress the  ${}^{\text {mot}}$ in  $\PP ^{\text {mot}}$
and simply write  $\PP $.

Recall that if   $Y$ is a smooth subvariety of the
smooth variety $X$
and $\widetilde  X$ denotes the blow-up of $X$ along $Y$ and  $E$ is the
exceptional divisor
$$
\PP \widetilde  X=
\PP X +\PP E-\PP Y
$$
In our case we get
$$\align
\PP \widetilde   M_j&=
\PP M_{j-1}+\PP E_j -
\PP \Bbb PW_j^-  \\
\PP \widetilde  M_j&=
\PP M_{j}+\PP E_j-
\PP \Bbb PW_j^+
\endalign $$
combining both equalities
$$
\PP M_j=\PP M_{j-1}+ \PP \Bbb PW_j^+-
\PP \Bbb PW_j^-
$$

Projective bundles are rationally cohomologically trivial so
$$\align
\PP M_j=&
\PP M_{j-1}+
\PP C^{(j)}(\1+ \cdots + \1(-d+2j-g+2))\\
&-
\PP C^{(j)}
(\1+ \cdots + \1(-j+1))
\endalign $$
this is
$$
\PP M_j=\PP M_{j-1}+\PP C^{(j)}
(\1(-j)+ \cdots + \1(-d+2j-g+2))
$$
when  $j=0$ this is still valid taking $M_{-1}=\emptyset $ since $M_0$
is
just $\Bbb P^{d+g-2}$. Now add all these expressions for $j=0$ to  $j=i$
$$
\PP M_i
=\underset j=0\to {\buildrel i\over \sum }
\PP C^{(j)}
(\1(-j)+ \cdots + \1(-d+2j-g+2))
$$
this proves the first part of the
proposition. For the rest re-write the last expression
$$
\PP M_i
=
\underset j=0\to {\buildrel i\over \sum }
\PP C^{(j)}
{\1(-j)-\1(-d+2j-g+1)\over \1-\1(-1)}
$$
and apply the corollary to our motivic MacDonald formula
$$
\align
\PP M_i&=
\SSS
\underset  T^{j}\to {\text {Coef}}
{(1+T)^{h^1C}\over (1-\1 T)(1-\1(-1)T)}
{\1(-j)-\1(-d+2j-g+1)\over \1-\1(-1)}   \\
&=
\underset T^{i}\to {\text {Coef}}
\SSS
{(1+T)^{h^1C}T^{i-j}  \over (1-\1 T)(1-\1(-1)T)}
{\1(-j)-\1(-d+2j-g+1)\over \1-\1(-1)}   \\
&=
\underset T^{i}\to {\text {Coef}}
\SSS
T^{i-j}(-j)-T^{i-j}(-d+2j-g+1)
\\
&\qquad \cdot
{\1\over \1-\1(-1)}
{(1+T)^{h^1C}\over (1-\1 T)(1-\1(-1)T)} \\
&=
\underset T^{i}\to {\text {Coef}}
\left(
{T^{i+1}-\1(-i-1)\over T-\1(-1)}+
{(\1 -T^{i+1}(-2i+2))\1 (-d-g+1+2i)\over \1(-2) T-\1}
\right)
\\
&\qquad
\cdot {\1\over \1-\1(-1)}
{(1+T)^{h^1C}\over (1-\1 T)(1-\1 (-1)T)}
\endalign
$$
the Proposition is
proved.
\hfill \qed \enddemo

In the odd degree case, if  $d>4g-4$  $M_\omega $ is a  $\Bbb
P^{d-2g+1}$-fibration over $N_0(2,d)$. As  $N_0(2,d)\simeq N_0(2,1)$ we
can choose any convenient value of $d$, we use  $d=4g-3$, then  $\omega
=2g-2$. Then as projective fibrations split in rational motives,

$$
\PP M_\omega =
{\1-\1(-2g+1)\over \1-\1(-1)}
\PP N_0.
$$

If we put the formula for  $\PP M_i$ into the previous expression we
obtain  $\PP N_0(2,1)$ in the form

$$ \align
{-\1(-g)\over \1-\1(-2g+1)}&
\underset T^{2g-2}\to {Coef}
{(1+T)^{h^1C}\over (\1-\1(-2)T)(\1-\1T)(\1-\1(-1)T) }
\\
&
+{\1(-2g+2)\over \1-\1(-2g+1)}
\underset T^{2g-2}\to {Coef}
{(1+T)^{h^1C}\over (\1-\1(2)T)(\1-\1T)(\1-\1(-1)T) }
\endalign $$

\noindent
Now our problem is to simplify this in  $K_0\MD $. We shall need a
definition.

\proclaim{Definition}
Let  $A,B$ and $M$ be objects in  $\Cal M^{AH}_k$ with
$r=\text{rank}M$,
define $(A+B)^M$ to be the Newton binomial
$$(A+B)^M=\sum \wedge ^iM \cdot A^{ r-i} \cdot
B^{ i}\in  K_0\MD  $$
\endproclaim

\demo{Caution} It is not true that $(A+B)^M=(B+A)^M$,
maybe it would be a better idea to write  $\lambda _{A,B}(M)$ instead of
$(A+B)^M$. However
in the case of interest  $M=h^1C$ there is a certain relation.
\enddemo

\proclaim{Main Lemma} If  $M=h^1C$,  with  $C$ a curve of genus  $g$, we
have $$(A+B)^M=(B(-1)+A)^M(g).$$
\endproclaim

\demo{Proof}
Poincar\'e duality on the Jacobian of  $C$ says  $\wedge ^iM\simeq (\wedge
^{2g-i}M(g))^{\vee }$, and by Poincar\'e duality on  $C$,  $M^{\vee
}\simeq M(1)$ so that
$$\align
\wedge ^iM&\simeq (\wedge ^{2g-i}M(g))^{\vee }\simeq (\wedge
^{2g-i}M^{\vee})(-g) \\
&\simeq  \wedge
^{2g-i}(M(1))(-g)\simeq \wedge ^{2g-i}M(g-i)
\endalign $$
apply this to the definition of  $(A+B)^M$
$$\align
(A+B)^M&=\wedge ^0MA^{2g}+\wedge ^1MA^{2g-1}B+\wedge
^2MA^{2g-2}B^2+\cdots +\wedge ^{2g}MB^{2g} \\
&=\wedge ^{2g}M(g)A^{2g}+\wedge ^{2g-1}M(g-1)A^{2g-1}B+\cdots
+\wedge ^0M(-g)B^{2g} \\
&=(\wedge ^{2g}MA^{2g}+\wedge ^{2g-1}M(-1)A^{2g-1}B+\cdots
+\wedge ^0M(-2g)B^{2g})(g)\\
&=(B(-1)+A)^M(g) .
\endalign $$
\hfill \qed
\enddemo

\proclaim{Theorem}
If  $N_0(2,1)$ denotes the moduli space of rank two
vector bundles with fixed odd degree on a curve  $C$ then its
motivic Poincar\'e polynomial in  $K_0(\MD )$ is
$$
\PP ^{\text {mot}} N_0(2,1)=
{
(\1+\1(-1))^{h^1C}-(\1
+\1)^{h^1C}(-g) \over
(\1-\1(-1))(\1-\1(-2))
}.
$$
\endproclaim

\demo{Proof of theorem}
We have seen that
$$
\align
\PP N_0&=
{-\1(-g)\over \1-\1(-2g+1)}
F(\1,\1(-1),\1(-2))            \\
&+
{\1(-2g+2)\over \1-\1(-2g+1)}
F(\1 ,\1(-1),\1(1))
\endalign
$$
where in analogy with \cite {Th}, $F(a,b,c)$ means
$$
F(a,b,c)=\underset T^{2g}\to {Coef}
{(\1+T)^{h^1C}\over (\1 -aT)(\1-bT)(\1-cT)}
$$
$a,b,c$ are now motives.
By direct calculation
one can prove the same identity as in
\cite {Th},
$$
F(a,b,c)=
{(a+\1)^{h^1C}\over (a-b)(a-c)}
+
{(b+\1)^{h^1C}\over (b-c)(b-a)}
+
{(c+\1)^{h^1C}\over (c-a)(c-b)}
$$
then  $\PP N_0$ equals

\centerline {$
{\2(-2g+2)\over \2-\2(-2g+1)}
\left(
{ (\2+\2)^{h^1C}\over (\2-\2(-1))(\2-\2(1))}
+
{ (\2(-1)+\2)^{h^1C}\over (\2(-1)-\2)(\2(-1)-\2(1))}
+
{ (\2(1)+\2)^{h^1C}\over (\2(1)-\2)(\2(1)-\2(-1))}
\right)
$}

\centerline {$
+{-\2(-g)\over \2-\2(-2g+1)}
\left(
{ (\2+\2)^{h^1C}\over (\2-\2(-1))(\2-\2(-2))}
+
{ (\2(-1)+\2)^{h^1C}\over (\2(-1)-\2)(\2(-1)-\2(-2))}
+
{ (\2(-2)+\2)^{h^1C}\over (\2(-2)-\2)(\2(-2)-\2(-1))}
\right)
$}

\medbreak

Call  $S_1$ the result of adding the third summand in both sums
and $S_2$ the rest, we shall first calculate  $S_1$,
$$
\align
{\2(-2g+2)(\2(1)+\2)^{h^1C}
\over
(\2(1)-\2)(\2(1)-\2(-1))}
&=
{\2(-2g+2)\2(2g)(\2+\2(-1))^{h^1C}\2(-1)\2(-1)
\over
(\2-\2(-1))(\2-\2(-2))} \\
&=
{(\2+\2(-1))^{h^1C}\over (\2-\2(-1))(\2-\2(-2)) }
\endalign
$$
and
$$
{-\2(-g)(\2(-2)+\2)^{h^1C}
\over (\2(-2)-\2)(\2(-2)-\2(-1)) }
=
{-\2(-2g)(\2+\2(-1))^{h^1C}\2(1)\over
(\2-\2(-1))(\2-\2(-2))}
$$

Adding and dividing by  $(\1-\2(-2g+1))$
$$
S_1={(\2+\2(-1))^{h^1C}\over
(\2-\2(-1))(\2-\2(-2))}
$$
similarly we calculate  $S_2$
$$
S_2=-{(\2+\2)^{h^1C}\2 (-g)\over
(\2-\2(-1))(\2-\2(-2))}
$$
sum  $S_1$ and  $S_2$ to get the desired expression for  $\PP N_0$.
This proves the theorem.\hfill \qed
\enddemo

The formula in the theorem contains the one found by Harder and
Narasimhan in \cite {HN},
$$
P_tN_0(2,1)={ (1+t^3)^{2g}-t^{2g}(1+t)^{2g}\over (1-t^2)(1-t^4) }
$$
but now we can also get the Hodge numbers and the level of the Hodge
structure.
Recall that the level of a Hodge structure is
$ \underset h^{p,q}\neq 0\to {Max} \vert p-q\vert $.

\proclaim{Corollary} The Poincar\'e-Hodge polynomial of  $N_0(2,1)$ is
$$
P_{xy}N_0(2,1)=
{
(1+x^2y)^g(1+xy^2)^g-x^gy^g(1+x)^g(1+y)^g
\over
(1-xy)(1-x^2y^2)
}
$$
and the level of the Hodge structure  $H^iN_0(2,1)$ is less than or
equal to $\left[ {i\over 3}\right] $. \endproclaim

\demo{Proof}
Let  $\bold {BiGrad\text{-}Vec}_{\Bbb C}$ or  $\bold {Vec}_{\Bbb
C}[x,y]$ denote the category of vector spaces with a double graduation
$$
V=\underset n\to \oplus \underset i+j=n\to \oplus V^{i,j}
$$

Sometimes we shall write  $V^{i,j}x^iy^j$ instead of  $V^{i,j}$ to
remind us of the graduation. Note that  $K_0\bold
{BiGrad\text{-}Vec}_{\Bbb C}=\Bbb Z [x,y,x^{-1},y^{-1}]$ and  $[V]=\sum
\text{dim}V^{i,j}x^iy^j$.

We have to study the image of
$$
{
(\1+\1(-1))^{h^1C}-\1(-g)(\1
+\1)^{h^1C} \over
(\1-\1(-1))(\1-\1(-2))
}
$$
by the morphism
$$
K_0\MD \longrightarrow K_0\bold{BiGrad\text{-}Vec}_{\Bbb C}=\Bbb
Z[x,y,x^{-1},y^{-1}] $$
as this is a morphism of rings it is enough to calculate the
image of ${\2\over \2-\2(-1)}$,
${\2 \over \2-\2(-2)}$,
$\1(-g)$,
$(\1+\1)^{h^1C}$ and
$(\1+\1(-1))^{h^1C}$. Consider the functor
$$
\MD \longrightarrow \bold{BiGrad\text{-}Vec}_{\Bbb C}
$$

It takes  $\1(-i)$ to  $\Bbb Cx^iy^i$ so the image of
${\2 \over \2-\2(-1)}$,
${\2 \over \2-\2(-2)}$ and
$\1(-g)$ by the morphism of  $K_0$ rings is just  ${1\over 1-xy}$,
 ${1\over 1-x^2y^2}$ and  $x^gy^g$.

This functor sends  $(\1+\1)^{h^1C}=\underset n\to \oplus \wedge ^nh^1C$
to  $\underset n\to \oplus \wedge ^n(\Bbb C^gx\oplus \Bbb C^gy)=
\underset n\to \oplus \underset i+j=n\to \oplus \wedge ^i(\Bbb
C^gx)\otimes \wedge ^j(\Bbb C^gy)$ and going down to
$K_0$ we get  $(1+x)^g(1+y)^g$.

Similarly we get the Poincar\'e-Hodge polynomial of the motive
$(\1+\1(-1))^{h^1C}$
and putting it all together we obtain the Poincar\'e-Hodge polynomial of
the moduli space.

Note that if  $A=(1+xy^2)(1+x^2y)$ and  $B=xy(1+x)(1+y)$ then
$$ \align
P_{xy}N_0(2,1)
&={A^g-B^g\over A-B}      \\
&=A^{g-1}+A^{g-2}B+\cdots +B^{g-1}
\endalign $$
and as the only monomials in  $A$ and  $B$ are  $x^iy^j$ with  $i=j$,
$i=2j$ or  $2i=j$ one can now see that the level of  $H^i$ is
less than or equal to $\left[ {i\over 3}\right] $
\hfill \qed
\enddemo

\bigbreak
\define\dm{{d\over 2}}

\heading
\S 5. The mixed motive of  $N_0(2,0)$
\endheading

In the even determinant case the moduli space  $N_0(2,d)$  is not smooth
and so the motive of  $N_0(2,d)$ is a mixed motive, in this section we
shall not find this motive but rather the pure motivic Poincar\'e
polynomial.

Recall that we have $M_\omega \longrightarrow N_0(2,d)$ and
the motive of  $M_\omega $ is known. If  $d$ is odd this is a  $\Bbb
P$-fibration whereas if  $d$ is even it is only a  $\Bbb P$-fibration
over the stable locus (=nonsingular locus if $g>2$), $N_0(2,d)^s$. Call
$M_\omega ^s$ its preimage
$$ \CD
M_\omega ^s   @>>>    M_\omega   \\
@VVV                @VVV          \\
N_0(2,d)^s    @>>>    N_0(2,d)
\endCD $$
The strictly semistable locus of  $N_0(2,d)$ is known to be
isomorphic to the Kummer variety of  $C$,  $Kum(C)={Jac(C)\over
x\approx
-x}$ (see \cite {BS1}), this is again a singular variety with  $2^{2g}$
ordinary double points.

So to calculate the pure motivic Poincar\'e polynomials the strategy
will be:

\item {} {\it Step 1.}
Calculate the motive of  $M_\omega -M_\omega ^s$.
\item {} {\it Step 2.}
To get then the motive of  $M_\omega ^s$
\item {} {\it Step 3.}
Use  $M_\omega ^s\longrightarrow N_0(2,d)^s$ is a  $\Bbb P$-bundle to
obtain motive of  $N_0(2,d)^s$
\item {} {\it Step 4.}
Combine  {\it Step 3.} with  $P_tKum(C)$ to calculate  $P_tN_0(2,d)$.

\medbreak
\subheading {{\it Step 1.} Motive of  $\bold{M_\omega -M_\omega ^s}$}

We are working with vector bundles with fixed even determinant, fix this
to be  $\Cal O(dP)$,  $d$ is the degree and  $P$ is a point of  $C$.
We are concerned with $$\align M_\omega -M_\omega ^s
&=\{ (E,s)\vert
\sigma \text{-stable pair with } E \text { strictly semistable} \}
\\
&=\left\{
(E,s)\vert                    {
{ E\text { is an extension }
0\longrightarrow
\Cal L(\dm P)
\longrightarrow E\longrightarrow
\Cal L^{-1}(\dm P)\longrightarrow 0}
\atop
{\text {with } s\notin H^0(C,\Cal L(\dm P)),\quad \text{deg} \Cal
L=0}                            }
\right\}
\endalign $$

Let us now code this information:
$$ \matrix
H^0E &\longrightarrow &H^0\Cal L^{-1}(\dm P) \\
s    &\longrightarrow & \gamma \neq 0
\endmatrix $$
$\gamma $ determines a divisor  $D$, linearly equivalent to  $\Cal
L^{-1}(\dm P)$, conversely any degree  $d$ divisor  $D$, gives a
$\Cal L=\Cal O(\dm P-D)$, this doesn't suffice to recover
$(E,s)$. Consider
$$ \matrix
0&\longrightarrow &\Cal L(\dm P)_{\vert D}
&\longrightarrow &E_{\vert D}
&\longrightarrow &\Cal L^{-1}(\dm P)_{\vert D}&\longrightarrow &0 \\
0&\longrightarrow &H^0\Cal L(\dm P)_{\vert D}
&\longrightarrow &H^0E_{\vert D}
&\longrightarrow &H^0\Cal L^{-1}(\dm P)_{\vert D}&\longrightarrow &0 \\
&&&& s_{\vert D} &\longrightarrow & \gamma _{\vert D}=0
\endmatrix $$
so we get a  $p\in H^0\Cal L(\dm P)=H^0\Cal O(dP-D)$.

\proclaim {Proposition}
 $D\in Div^\dm (C)=C^{(\dm )}$ and  $p\in \Bbb PH^0\Cal O(dP-D)$
determine  $(E,s)$ uniquely.
\endproclaim

\demo {Proof} See \cite {Th, (3.3)}
\enddemo

\demo {Remark} These data are parametrized by  $\Bbb PW^-_\dm $, a
certain  $\Bbb P^{\dm -1}$-bundle over  $C^{(\dm )}$.
\enddemo

\proclaim {Corollary} The inverse image of the singular locus by
 $M_\omega \longrightarrow N_0(2,\Cal O(dP))$ is isomorphic to a  $\Bbb
P^{\dm -1}$-bundle over  $C^{(\dm )}$.
\endproclaim

\demo {Remarks}
{\it 1.} Recall from \S 3  that
$$\align
\text {dim}N_0=3g-3 \hskip 3truecm & \text {dim}M_\omega
-M_\omega^s=\dm -1+\dm=d-1 \\
\text {dim}M_\omega =d+g-2 \hskip 3truecm & \text
{codim}_{M_\omega } (M_\omega -M_\omega^s)=g-1
\endalign $$

{\it 2.} $M_\omega \longrightarrow N_0(2,d)$ is surjective for
$d>2g-2$ for then every semistable bundle has nonzero sections. Also
$H^1E=0$ if  $E$ stable for  $d>4g-4$. So for  $d$ odd one uses
$d=4g-3$ and then  $\omega =[{d-1\over 2}]=2g-2$, for  $d$ even we use
$d=4g-2$ and  $\omega =[{4g-3\over 2}]=2g-2$.
\enddemo

We saw that  $M_\omega -M_\omega ^s$ is isomorphic to a  $\Bbb
P^{\dm-1}$-bundle over  $C^{(\dm)}$ so
$$ P_t(M_\omega -M_\omega ^s)=
P_t\Bbb P^{\dm-1}P_tC^{(\dm)}=
P_t\Bbb P^{2g-2}P_tC^{(2g-1)}
$$
now we don't need any MacDonald formula as  $C^{({2g-1})}$ is a  $\Bbb
P^{g-1}$-bundle over  $Jac(C)$ so
$$ P_t(M_\omega -M_\omega ^s)=
{\1-\1(1-2g)\over \1-\1(-1)}
{\1-\1(-g)\over
\1-\1(-1)}(\1+\1)^{h^1C}$$

\medbreak
\subheading {{\it Step 2.} Motive of  $\bold {M_{\omega }^s}$}

Recall that $\PP M_i$ equals
$$
{\1\over \1-\1(-1)}
\underset T^{i}\to {\text {Coef}}
\left(
{\1(-d+2i-g+1)\over \1(-2)T-\1}-
{\1(-i-1)\over T-\1(-1)}
\right)
{(1+T)^{h^1C}\over (\1-\1 T)(\1-\1(-1)T)}
$$
In our case  $i=\omega =2g-2$ and  $d=4g-2$ so that  $-d+2i+1-g=-g-1$
and  $-i-1=-2g+1$, so
$$\align
P_t M_\omega &=
{\1\over \1-\1(-1)}
\underset T^{2g-2}\to {\text {Coef}}
\left(
{\1(-g-1)\over \1(-2)T-\1}-
{\1(-2g+1)\over T-\1(-1)}
\right)
{(1+T)^{h^1C}\over (\1-T)(\1-\1(-1)T)} \\
&={-\1(-g-1)\over \1-\1(-1)}F(\1
,\1(-1),\1(-2))+{\1(2-2g)\over \1-\1(-1)}F(\1
,\1(-1),\1(1))\\
&={\1(-g)-\1(-g-1)\over \1-\1(-1)}F(\1
,\1(-1),\1(-2))+\PP M_\omega (2,1) \endalign $$
where  $M_\omega (2,1)$ indicates the moduli space for odd degree
considered in the previous section, we know its Poincar\'e polynomial
$$\align
P_tM_\omega (2,1)&=
{\1-\1(-2g+1)\over \1-\1(-1)}
P_tN_0(2,1)\\
&=
{\1-\1(-2g+1)\over \1-\1(-1)}
{ (\1+\1(-1))^{h^1C}-\1(-g)(\1+\1)^{h^1C} \over (\1-\1(-1))(\1-\1(-2)) }
\endalign $$
Now  $P_tM_\omega $ is
$$ \1(-g)F(\1 ,\1(-1),\1(-2))+
({\1-\1(-2g+1)})
{ (\1+\1(-1))^{h^1C}-\1(-g)(\1+\1)^{h^1C} \over
(\1-\1(-1))^2(\1-\1(-2)) } $$
and we use the formula for  $F(a,b,c)$ to obtain
$$\align
{}
&{}_{
\1(-g)
\left(
{ (\1+\1)^{h^1C}\over (\1-\1(-1))(\1-\1(-2)) } +
{ (\1(-1)+\1)^{h^1C}\over (\1(-1)-\1)(\1(-1)-\1(-2)) } +
{ (\1(-2)+\1)^{h^1C}\over (\1(-2)-\1)(\1(-2)-\1(-1)) }
\right)
}
\\
&{}_{
+
(\1-\1(-2g+1))
{ (\1+\1(-1))^{h^1C}-\1(-g)(\1+\1)^{h^1C} \over (\1-\1(-1))^2(\1-\1(-2))
}  }
\endalign
$$

Doing a bit of calculation we get

\medbreak
\centerline {$
\PP M_\omega=
{(\1+\1(-1))^{h^1C}(\1-\1(-2g))-
\1(-g-1)(\1+\1)^{h^1C}(\1+\1(-g))(\1+\1(2-g))
\over
(\1-\1(-1))^2(\1-\1(-2))}
$}
\medbreak

Now is when pure Poincar\'e polynomials of mixed motives enter the
scene,

\medbreak
\noindent
$\PP M_\omega ^s
=P_tM_\omega -\1(1-g)\cdot P_t(M_\omega -M_\omega ^s) $

\noindent
\hskip 1.2truecm $=
{
(\1+\1(-1))^{h^1C}(\1-\1(-2g))-\1(-g-1)(\1+\1)^{h^1C}
(\1+\1(-g))(\1+\1(2-g))
\over
{}            }$

\noindent
\hskip 1.7truecm ${
-\1(1-g)(\1+\1)^{h^1C}(\1-\1(-g))(\1-\1(1-2g))(\1-\1(-2))
\over
(\1-\1(-1))^2(\1-\1(-2)) }
$

\medbreak
\subheading {{\it Step 3.} Motive of  $\bold {N_0(2,0)^s}$}

As  $M_\omega ^s\longrightarrow N_0(2,0)^s$ is a  $\Bbb
P^{2g-1}$-bundle, the cohomology of the  $N_0(2,0)^s$ will be that of
$M_\omega ^s$ divided by  ${\2-\2(-2g)\over \2-\2(-1)}$, that is
$\PP N_0(2,d)^s$ equals
$$
{}_{
{ (\2+\2(-1))^{h^1C}\over
(\2-\2(-1))(\2-\2(-2))     }
-
{
\2(-g)(\2+\2)^{h^1C}
\left(
\2(-1)(\2+\2(-g))(\2+\2(2-g))
+\2(1)(\2-\2(-g))(\2-\2(1-2g))(\2-\2(-2))
\right)
\over
(\2-\2(-1))(\2-\2(-2))(\2-\2(-2g))
}
}
$$

\medbreak
\subheading {{\it Step 4.} Motive of  $\bold {N_0(2,0)}$}

Use the previous expression together with
$$
\PP N_0(2,0)=\PP N_0(2,0)^s+\PP Kum(C) (-2g+3)
$$
Note that the motivic Poincar\'e polynomial of the Kummer variety of
$C$ is given by
$$P_tKum(C)={1\over 2}\left( (\1+\1)^{h^1C}+(\1-\1)^{h^1C}\right)
=\underset i=0\to {\buildrel g\over \oplus } \wedge^{2i}h^1C
.$$

\proclaim{Corollary}
The mixed AH motive  $h^iN_0(2,0)$ has only weights  $i$ and  $i+1$.
Furthermore, if $i<2g-2$, $h^iN_0(2,0)\simeq h^iN_0(2,1)$, in
particular it is pure of weight $i$. \endproclaim

\demo{Proof}
$M_\omega $ and  $M_\omega -M^s_\omega $ are smooth and
projective so their cohomology is pure, now by Gysin
$$
{}_{
\cdots\longrightarrow
h^{i-2g+2}(M_\omega -M_\omega
^s)(-g+1)
\longrightarrow
h^iM_\omega
\longrightarrow
h^iM_\omega ^s
\longrightarrow
h^{i+1-2g+2}(M_\omega
-M_\omega ^s)(-g+1)
\longrightarrow \cdots
}
$$
so  $h^i(M_\omega ^s)$ has weights  $i$ and  $i+1$. As  this is a
projective bundle over $N_0(2,0)^s$, the same
happens with  $N_0(2,0)^s$ and the Gysin exact sequence asociated to
$Kum(C)\hookrightarrow N_0(2,0)$
$$
{}_{
\cdots\longrightarrow
h^{i-2g+6}Kum(C)(-2g+3)
\longrightarrow
h^i N_0(2,0)
\longrightarrow
h^i N_0(2,0) ^s
\longrightarrow
h^{i+1-2g+6}Kum(C)(-2g+3)
\longrightarrow \cdots
}
$$
gives the result.

The ring  $K_0 \MD $ is a graded ring so it makes sense to take the
$n$-truncation of an element of  $K_0 \MD $, going down to Poincar\'e
polynomials this is just taking the residue modulo  $t^n$.

The  $2g-2$ truncation of  $\PP N(2,0)$ coincides with the  $2g-2$
truncation of
$$
{ (\1 +\1(-1))^{h^1C} \over
(\1-\1(-1))(\1-\1(-2)) }
$$
and this coincides with that of  $\PP N(2,1)$.
\hfill \qed
\enddemo

\demo{Remark} Of course if the genus of the curve is 2 then the moduli
space  $N_0(2,0)$ is smooth so one gets the pure motive of this.
\enddemo

\bigbreak
\heading
\S 6. Jacobian considerations.
\endheading

Let  $X$ be a smooth projective complex variety. Griffiths
defines the $i$-th intermediate jacobian of $X$ to be:
$$
J^i(X)={
H^{2i-1}(X,\Bbb C)
\over
F^iH^{2i-1}(X,\Bbb C)+H^{2i-1}(X,\Bbb Z)
}
$$
This is just a complex torus with a naturally defined 2-form, it is not
an abelian variety in general.

Now let  $k$ be a field as before. If  $X$ is a smooth projective
variety over  $k$ for each embedding  $\sigma :k\hookrightarrow \Bbb C$
we get a complex variety  $X_\sigma $ and intermediate jacobians for
each  $X_\sigma $.

Note that if we know $h^{2i-1}(X)$ and
are interested in  $J^i(X_{\sigma })$ there is only one piece of data
missing: the
entire structure on the singular cohomology group
$H^{2i-1}(X_\sigma ,\Bbb
Q)$, so we can recover  $J^i(X_\sigma )$  up to isogeny from
$h^{2i-1}(X)$

\medbreak

As we have just seen one can recover the intermediate
jacobians up to isogeny from the motive of the variety. We shall now
exploit the expression found for the  motivic Poincar\'e polynomial of
the moduli spaces  $N_0(2,1)$ and  $N_0(2,0)$ to get information on
the intermediate jacobians.

\proclaim{Corollary} Let  $\delta \in \{0,1\} $. For  $i\leq g-1+\delta
$ the $i$-th intermediate jacobian of the moduli
space $N_0(2,\delta )$  is isogenous to
$$
\underset \alpha =1\to {\buildrel \left[ {i+1\over 3}\right] \over
\prod } J^{\alpha }Jac(C)^{\left[ {i+3-2\alpha \over 2}\right] }
$$
\endproclaim

\demo{Proof}
We have proved that  $h^iN_0(2,0)\simeq h^iN_0(2,0)$
if  $i<2g-2$ so that it suffices to prove the theorem for the odd degree
case

In the smooth case, expand the formula in the theorem in power series to
get $$
((\1 +\1(-1))^{h^1C}-\1(-g)(\1 +\1)^{h^1C})
(\sum \1(-i))
(\sum \1(-2i))
$$
taking the  $2g$-truncation
$$
\underset n\geq 0\to \sum \wedge ^ih^1C(-i)
(\1+\1(-1)+\1(-2)+\cdots )
(\1+\1(-2)+\1(-4)+\cdots )
$$

The last expression is seen to be equal to
$$
\underset n\geq 0\to \sum \wedge ^ih^1C(-i)
\underset n\geq 0\to \sum \left[ {n\over 2}+1\right] \1(-n)
$$

Multiplying the series we get a weight decomposition
$$
\underset n \to \sum
\left(
\underset 3a+2b=n\to \sum
\wedge ^ah^1C(-a)\otimes \left[ {b\over 2}+1\right] \1(-b)
\right)
$$

As we are interested in odd cohomology  $n=2i-1$, and if  $3a+2b=2i-1$
then  $a$ has to be odd,  $a=2\alpha -1$  $\alpha \geq 1$. And the
condition  $3a+2b=n$ turns  into $3\alpha +b=i+1$. So that for  $i\leq
g$ we get
$$\align
h^{2i-1}N_0&
=\underset {\alpha =1} \to
{    \buildrel {\left[ {i+1\over 3}\right] } \over \sum   }
  \wedge ^{2\alpha -1}   h^1C(1-2\alpha )
\otimes \left[ {i+1-3\alpha \over 2} +1 \right] \1 (3\alpha -i-1) \\
&=\underset {\alpha =1} \to
{ \buildrel {\left[ {i+1\over 3}\right] } \over \sum }
\wedge ^{2\alpha -1}  h^1C(\alpha -i)^{\oplus \left[{i+3-3\alpha \over
2}\right] }
\endalign $$
now the result follows.\hfill \qed
\enddemo

\demo{Remarks}

-Newstead and Mumford proved in \cite {MN} that
$J^2N_0(2,1)$ is in fact isomorphic to  $Jac(C)$ as polarized abelian
varieties, as a corollary they obtain a Torelli-type theorem:  if
$N_0(2,1)$ is isomorphic for two curves then the curves are isomorphic.

-Balaji has calculated the second intermediate jacobian of the canonical
desingularization of  $N_0(2,0)$, $M$ (\cite {Ba}), there is a
canonical isogeny of degree  $2^{2g}$  $Jac(C)\longrightarrow J^2M$.
As in the case treated by Mumford and Newstead he obtains
a Torelly-type theorem.

-Our corollary applies as well to the $\ell $-adic intermediate
jacobians.
For the case  $i=2$ treated by Mumford and Newstead one can
re-do their proof in the  $\ell $-adic setup to obtain an isomorphism
and not only an isogeny.
\enddemo

\bigbreak

\enddocument